\documentclass{amsart}
\usepackage{latexsym}
\usepackage{amsmath}
\usepackage{amssymb}
\usepackage{amsfonts}
\usepackage{amsthm}
\usepackage{amsrefs}

\newcommand{\namedRef}[1]{%
\csname env:#1\endcsname\ \ref{#1}%
}

\def\envLRT#1#2{\expandafter\gdef\csname env:#1\endcsname{#2}}

\makeatletter
\newcommand{\namedLabel}[1]{\expandafter\xdef\csname env:#1\endcsname{\envName}%
\immediate\write\@auxout{\string\envLRT{#1}{\envName}}%
\label{#1}}
\newcommand{\NamedLabel}[1]{\expandafter\xdef\csname env:#1\endcsname{\envName}%
\immediate\write\@auxout{\string\envLRT{#1}{\envName}}%
\label{#1}}
\makeatother

\def\newType#1{\newtheorem{#1}[thm]{#1\gdef\envName{#1}}}
\newtheorem{thm}{Theorem}[section]
\newtheorem{Theorem}[thm]{Theorem\gdef\envName{Theorem}}
\newType{Lemma}

\theoremstyle{definition}
\newtheorem{Remark}[thm]{Remark\gdef\envName{Remark}}

\newcommand{\rightlabeledarrow}[3][]{%
\def\LRTxx{#1}\ifx\LRTxx\empty\relax%
\setbox0=\hbox{$\mathop{\expandafter\hbox%
{\rightarrowfill}}\limits_{\hbox{\ $\scriptstyle#3$\ }}^{\hbox{\ $\scriptstyle#2$\ }}$}%
\dimen0=\wd0\advance\dimen0 by 4pt%
\edef\LRTxx{\the\dimen0}\else%
\setbox0=\hbox to \LRTxx{$\mathop{\expandafter\hbox%
{\rightarrowfill}}\limits_{\hbox{$\scriptstyle#3$}}^{\hbox{$\scriptstyle#2$}}$}%
\dimen0=\wd0\advance\dimen0 by 4pt%
\fi%
\setbox0=\hbox{$\mathop{\expandafter\hbox%
{\rightarrowfill}}\limits_{\hbox{$\scriptstyle#3$}}^{\hbox{$\scriptstyle#2$}}$}%
\ifdim\dimen0<\wd0 \relax\dimen0=\wd0 \advance\dimen0 by 4pt\fi%
\ \hbox to\dimen0{\hfill\hbox to 0pt{\hss$\mathop{\hbox to\LRTxx {\rightarrowfill}}%
\limits_{\hbox to 0pt{\hss{$\scriptstyle#3$}\hss}}^{\hbox to 0pt{\hss{$\scriptstyle#2$}\hss}}$\hss}\hfill}\ %
}

\def\disjointunion{\perp\penalty10000\hskip-5pt\perp}
\def\Z{\mathbb Z}
\def\cy#1{\Z/{#1}\Z}

\def\pin{\ifmmode Pin^{-}\else$Pin^{-}$\fi}
\def\spin{\ifmmode Spin\else\hbox{$Spin$}\fi}
\def\qe{q}
\def\modTwo#1#2{#1 \bullet #2}
\def\standinc#1{2\cdot{#1}}
\def\MSpin{\text{\bf MSpin}}

\begin{document}
\title{Quadratic enhancements of surfaces: two vanishing results}
\author{L.~R.~Taylor}
\address{Department of Mathematics
 \newline\indent
University of Notre Dame
Notre Dame, IN 46556, USA}
\email{taylor.2@nd.edu}

\begin{abstract}
This note records two results which were inexplicably omitted from 
our paper on Pin structures on low dimensional manifolds, \cite{KT}. 
Kirby chose not to be listed as a coauthor.

A \pin{--}structure on a surface $F$ induces a quadratic enhancement of the mod $2$ 
intersection form,
$\qe\colon H_1(F;\cy{2}) \to \cy{4}$.

Theorem 1.1 says that $\qe$ vanishes on the kernel of the map in homology 
to a bounding $3${--}manifold. 
This is used by Kreck and Puppe \cite{KP} who refer for a proof to an email of 
the author to Kreck.
A more polished and public proof seems desirable. 

In \cite{KT}*{Section 6}, a \pin{--}structure is constructed on a surface 
$F$ dual to $w_2$ in an oriented $4${--}manifold, $M^4$.
Theorem 2.1 says that $\qe$ vanishes on the Poincar\'e dual 
to the image of $H^1(M^4;\cy{2})$ in $H^1(F;\cy{2})$.
\end{abstract}

\maketitle

\section{Surfaces bounding $3${--}manifolds}
Recall that a  fixed \pin{--}structure on a surface $F$
defines a function
$$\qe\colon H_1(F;\cy{2}) \to \cy{4}$$
which satisfies the formula 
$$\qe(x+y) = \qe(x) + \qe(y) + \standinc{(\modTwo x y)}
\text{  for all $x$, $y \in H_1(F;\cy{2})$.}$$ 
Here $\modTwo x y$ denotes the mod{--}$2$ intersection of the two classes 
and $\standinc\colon \cy{2}\to\cy{4}$ denotes the standard inclusion.
Conversely, every such function comes from a unique \pin{--}structure. 
This is classical, but see \cite{KT}*{Section 3}.

\begin{Theorem}
Let $M^3$ be a $3${--}manifold with a fixed \pin{--}structure and 
let $F$ be the boundary of $M$. 
Give $F$ the induced \pin{--}structure. 
Let $x\in H_1(F;\cy{2})$ be a class which vanishes in $H_1(M;\cy{2})$.
Then $\qe(x) = 0$. 
\end{Theorem}

\begin{proof}
We start with two lemmas.
 
\begin{Lemma}\namedLabel{can do surgery}
Let $S^1\subset F$ be an embedded circle with trivial normal bundle.
Fix a \pin{--}structure on $F$. 
One can do surgery on this embedding and extend the \pin{--}structure 
to the trace of the surgery if and only if $\qe(S^1) = 0$.
\end{Lemma}
\begin{proof}
Pick a point on the circle and orient the tangent space at this point and also orient 
the circle at the point.
This orients, and hence frames, the normal bundle to $S^1$ in $F$. 

A tubular neighborhood of the circle is now oriented and so a \pin{--}structure on 
$F$ restricts to a \spin{--}structure on this neighborhood. 
The framing on the normal bundle induces a stable framing of $S^1$ and
$\qe(S^1) \in \MSpin_1 \cong\cy{2}$.

This is equivalent to the description in Kirby{--}Taylor \cite {KT} just before Definition 3.5.
The definition given  there works for all circles, not just ones with trivial 
normal bundle.

The trace of the surgery is formed by gluing $D^2\times D^1$ to $S^1\times D^1$. 
Since $1${--}dimensional framed bordism is $\cy{2}$ and maps isomorphically to $\MSpin_1$, 
if $\qe(S^1) = 0 \in \MSpin_1$ then the \spin{--}structure on the circle extends over a disk
and hence over $D^2\times D^1$ and finally over the entire trace.

The only if part follows from the next lemma.
\end{proof}

\begin{Lemma}\namedLabel{embedded disk}
Let $X$ be a surface bounding a $3${--}manifold $W$ and suppose $W$ has a \pin{--}structure.
Let $S^1\subset F$ be an embedded circle which bounds an embedded disk in $W$.
Then $\qe(S^1) = 0$.
\end{Lemma}
\begin{proof}
A tubular neighborhood of the disk in $W$ is trivial.
As in the proof of the first lemma, orient a point on the circle and orient the circle.
These orientations extend over the neighborhood of the disk and over the disk and 
hence frame the normal bundle of the disk in $W$. 
The \pin{--}structure on $W$ restricts to a \spin\ structure on the neighborhood of the disk.

Lemma 2.7 of \cite{KT} shows that restricting to the disk and then to the bounding circle 
gives the same \spin{--}structure as restricting to $X$ and then to the circle.
The first restriction is obviously $0$ and the second is $\qe(S^1)$.
\end{proof}

Turn now to the proof of the Theorem.
\begin{proof}
Assume $M$ is connected and hence has a handlebody decomposition with no $0${--}handles.
Divide $M$ into two pieces, $Y$ and $M^\prime$.
The submanifold $Y$ is obtained from $F$ by attaching the $1${--}handles. 
Let $X$ be the rest of the boundary of $Y$.
Let $M^\prime$ be the result of attaching the $2$ and $3$ handles to $X$ so that 
$M = Y \cup M^\prime$. 
Notice that $Y$ is obtained from $X$ by adding $2${--}handles. 
The \pin{--}structure on $M$ restricts to one on $Y$ and one on $M^\prime$ and hence
one on $X$. 
Let $\qe_X$ denote the resulting quadratic enhancement.

Let $x\in H_1(F;\cy{2})$ be a class that vanishes in $H_1(M;\cy{2})$.  
The class $x$ can always be represented by disjoint embedded circles using the usual
trick for removing transverse intersections. 
Let $\kappa\in H_2(M, F;\cy{2})$ be the resulting relative class. 
By excision, $H_2(M,F;\cy{2}) = H_2(M^\prime,X;\cy{2})$ so let
$x_1\in H_1(X;\cy{2})$ be the boundary of $\kappa$ after excision. 
It follows from the construction that $x_1$ vanishes in $H_1(M^\prime;\cy{2})$ and
that there is a relative class $\lambda\in H_2(Y, F\disjointunion X;\cy{2})$ 
with boundary $x + x_1$.

Every homology class in $X$ which dies in $M^\prime$ can be represented 
after handle slides and additions by the boundary of a 2{--}handle. 
But this is a surgery so by \namedRef{can do surgery}, $\qe_X(x_1) = 0$.

By adding all the $1${--}handles in a small disk in $F$, we see $X$ as a connected sum 
of $F$ and some tori and Klein bottles.
In particular, there is a class $\bar{x}$ in $H_1(X;\cy{2})$ which can be joined to 
$x\in H_1(F;\cy{2})$ by an embedded cylinder in $Y$.
Then $\bar{x} + x_1$ vanishes in $H_1(Y;\cy{2})$ and again $Y$ has no $0$ or $1${--}handles
when built from $X$, so $\qe_X(\bar{x} + x_1) = 0$.

Make the cylinder from $x$ to $\bar{x}$ transverse to a surface spanning $x$ and $x_1$ 
representing $\lambda$. 
The intersection will be some circles and some arcs and the usual ``arc has two ends'' 
argument shows 
$\modTwo{\bar{x}}{x_1} = \modTwo{x}{x}$. 
But since $x$ bounds in $M$, $\modTwo{x}{x} = 0$. 
Hence $0 = \qe_X(\bar{x} + x_1) = \qe_X(\bar{x}) + \qe_X(x_1) = \qe_X(\bar{x})$. 

Since we added the $1${--}handles in a small disk $D^2\subset F$, we see an embedding
$F^\prime\times [0,1] \subset Y$, where $F^\prime = F - D^2$. 
Furthermore, this embedding is the inclusion $F^\prime\subset F$ 
at $F^\prime \times 0$ and is the inclusion $F^\prime \subset X$ at $F^\prime\times 1$.
Here $F^\prime \subset X$ is the inclusion whose complement is 
the connected sum of tori and Klein bottles. 
The cylinder between our representatives of $x$ and $\bar{x}$ lies in $F^\prime \times [0,1]$. 

Since the \spin{--}structure on the circles representing $x$ or $\bar{x}$ can be computed 
from the \pin{--}structure on a neighborhood of these circles \cite{KT}*{Definition 3.5}, 
we can work in $F^\prime \times[0,1]$ with its induced \pin{--}structure. 

The \pin{--}structures on $F^\prime\times 0$ and $F^\prime\times 1$ are equivalent and
since equivalent \pin{--}structures induce the same restriction to an $S^1\subset F^\prime$,
$\qe_X(\bar{x}) = \qe(x)$ and hence $\qe(x)=0$.
\end{proof}

\section{The dual to $w_2$}
Let $M^4$ be an oriented $4${--}manifold and let $F\subset M$ be dual to $w_2$.
In \cite{KT}*{Section 6} we used a \spin{--}structure on $M - F$ to construct
a \pin{--}structure on $F$. 
\begin{Theorem}
The composition 
$$H^1(M;\cy{2}) \to H^1(F;\cy{2}) \rightlabeledarrow{\quad \cap [F]\quad }{}
H_1(F;\cy{2}) \rightlabeledarrow{\qe}{} \cy{4}$$ 
is trivial. 
\end{Theorem}
\begin{Remark}
This says that the image of $H^1(M;\cy{2})$ in $H^1(F;\cy{2})$ can be at most half{--}dimensional. 
\end{Remark}
\begin{proof}
The set of \spin{--}structures on $M-F$ is an $H^1(M;\cy{2})${--}torsor.
The set of \pin{--}structures on $F$ is an $H^1(F;\cy{2})${--}torsor.
By Lemma 6.2 of \cite{KT}, if the \spin{--}structure is changed by $x\in H^1(M;\cy{2})$,
then the \pin{--}structure changes by $i^\ast(x)$, where $i\colon F\to M$ denotes the embedding.

Associated to $\qe$ there is an element $\beta(\qe) \in \cy{8}$ called the Brown{--}Arf invariant. 
The only properties needed of this invariant will be recalled below.
Theorem 6.3 of \cite{KT} says that
$2\cdot \beta(F) = F \bullet F - \text{sign}(M)$ mod $16$, 
where $F\bullet F\in \Z$ is the self{--}intersection
and $\text{sign}(M)$ is the signature.  

Lemma 3.7 of \cite{KT} 
says that if the \pin{--}structure is changed by $y\in H^1(F;\cy{2})$,
then $\beta$ changes by $2\cdot \qe(y \cap [F])$. 

If the \spin{--}structure on $M-F$ is changed by $x\in H^1(M;\cy{2}$, then the right{--}hand side of
$2\cdot \beta(F) = F \bullet F - \text{sign}(M)$ does not change, hence neither can the left. 
Therefore $2\cdot \qe\bigl(i^\ast(x)\cap[F]\bigr) = 0$. 
But since $2\cdot\colon \cy{4}\to \cy{8}$ is injective, the result follows. 
\end{proof}

\end{proof}


\begin{thebibliography}{WW} 

\bib{KT}{article}{
   author={Kirby, R. C.},
   author={Taylor, L. R.},
   title={${\rm Pin}$ structures on low-dimensional manifolds},
   conference={
      title={Geometry of low-dimensional manifolds, 2},
      address={Durham},
      date={1989},
   },
   book={
      series={London Math. Soc. Lecture Note Ser.},
      volume={151},
      publisher={Cambridge Univ. Press},
      place={Cambridge},
   },
   date={1990},
   pages={177--242},
   review={\MR{1171915 (94b:57031)}},
}

\bib{KP}{article}{
  author = {Kreck , Matthias},
  author = {Puppe , Volker},
  title = {Involutions on 3-Manifolds and Self-dual, Binary Codes},
  eprint = {arXiv:0707.1599 [math.AT]},
  date = {2007}
}
\end{thebibliography}
\end{document}